# CONSISTENCY PROPERTIES OF A SIMULATION-BASED ESTIMATOR FOR DYNAMIC PROCESSES[1]

By Manuel S. Santos

*University of Miami and Universidad Carlos III de Madrid*

This paper considers a simulation-based estimator for a general class of Markovian processes and explores some strong consistency properties of the estimator. The estimation problem is defined over a continuum of invariant distributions indexed by a vector of parameters. A key step in the method of proof is to show the uniform convergence (a.s.) of a family of sample distributions over the domain of parameters. This uniform convergence holds under mild continuity and monotonicity conditions on the dynamic process. The estimator is applied to an asset pricing model with technology adoption. A challenge for this model is to generate the observed high volatility of stock markets along with the much lower volatility of other real economic aggregates.

**1. Introduction.** Simulation-based estimation is advocated in several applied areas of economics and finance (e.g., [6, 12, 15]), but relatively little is known about asymptotic properties of these estimation methods. Our purpose here is to establish some strong consistency properties of a simulation-based estimator. Although the estimator seems highly specific, our results should be of broad application to other types of simulation-based estimation methods. The estimator is applied to a macroeconomic model of technology adoption where some parameters are hard to estimate by other methods.[2] In this model, news about the arrival of new technologies will suddenly impact the stock market because of their option value in the creation of new products and designs. A challenge for this model is to reconcile the volatility of

---

Received September 2007; revised February 2009.

[1]Supported in part by Projects SEJ2005-05831 and ECO2008-04073 of the Spanish Ministry of Education, and a Cátedra de Excelencia of Banco Santander.

*AMS 2000 subject classifications.* Primary 62M05, 60K35; secondary 65C20, 60B10.

*Key words and phrases.* Markov process, simulation-based estimation, invariant probability, sample distribution, monotonicity, strong consistency.

[2]This is actually a typical situation, since models are abstractions of reality and so their parameters cannot usually be estimated from raw data.

real economic aggregates (e.g., worked hours, consumption and investment) with the much greater volatility of the stock market. Traditional business cycle models fail to offer a joint explanation for the fluctuations of the real and financial sectors. Hence, it becomes of interest to search for the best fit of our model and check if it has the ability to replicate volatilities along these two dimensions.

For simplicity, our theoretical analysis centers on the following parameterized family of random dynamical systems:

$$(1.1) \qquad s_n = \varphi(s_{n-1}, \varepsilon_n, \theta), \qquad n = 1, 2, \ldots.$$

Stochastic equations of the form (1.1) often arise as solutions of various dynamic models in biology, economics, physics and other applied disciplines. In the sequel, $s_n$ is a finite vector of state variables, $\varepsilon_n$ is a vector of i.i.d. shocks and $\theta$ is a finite vector of parameters. In many applications, it becomes crucial the estimation of the vector of parameters $\theta$ by a reliable procedure. Traditional data-based estimation methods are of limited applicability for nonlinear stochastic systems. First, the estimator may require optimization of the likelihood function, or of an associated objective, and such calculation may not be computationally feasible. Indeed, function $\varphi$ may be unknown, and usually has to be computed numerically. Besides, both the vector of shocks $\varepsilon$ and some state variables $s$ may be unobservable.[3] Second, commonly used goodness-of-fit criteria, such as the minimum size of the squared residuals, can be quite uninformative on the dynamic properties of the selected model. Simulation-based estimators can target directly the moments of an invariant probability or some other critical aspects of the dynamics.

The informational requirements for the implementation of our estimator are quite weak, and stem from the ability to simulate the model. Thus, a researcher must be able to evaluate function $\varphi$ (or have access to a computer-generated law of motion) as specified in (1.1), and take a stand on the functional form of the probability law for the shock $\varepsilon$. (Vector $\theta$ may encompass parameter values for the distribution of $\varepsilon$.) Therefore, the functional form of $\varphi$ and the actual sequence of realizations $\{\varepsilon_n\}_{n\geq 1}$ may both be unknown.

The proofs of consistency of the estimator deal with the nonlinear dynamical effects of the vector of parameters $\theta$, which feeds into the evolution of sample paths $\{s_n\}_{n\geq 1}$ for a fixed initial condition $s_0$ and a sequence of shocks $\{\varepsilon_n\}_{n\geq 1}$. Indeed, the estimation problem is defined over a continuum of invariant probability measures $\mu_\theta$ which vary with parameter $\theta$, and the estimated parameter values are selected over a continuum of sample distributions. In contrast, traditional estimators [16, 20] select values from the unique invariant distribution of the data generating process. A key step below is to demonstrate the uniform convergence of a continuum of sample distributions for every initial condition $s_0$ and almost all sequences of realizations $\{\varepsilon_n\}_{n\geq 1}$. This result is of independent interest in probability theory. As in usual proofs of consistency, the result amounts to a uniform law of large numbers, but for our simulation-based estimator the convergence is over a continuum of invariant distributions parameterized by vector $\theta$. This uniform law of large numbers is shown to hold under mild continuity and monotonicity conditions. We leave for future investigation extensions of this analysis to nonmonotone dynamical systems.

---

[3] These unobservable state variables are often called "latent variables." Note that function $\varphi$ could be nonseparable in $s$ and $\varepsilon$, and the ergodic sets and distributions of the stochastic dynamical system (1.1) will depend on parameter vector $\theta$.

The paper contributes to several strands of the literature. There are various related results on uniform convergence of invariant distributions for families of functions (e.g., [9, 25, 31]), but these results fall short of what is required in the present case since they are restricted to a fixed empirical distribution. Monotone Markov processes are studied in [2] and [10]. The simulation-based estimator is set forth in [11] and [21]. Available proofs on consistency of this estimator require uniform continuity and contractivity conditions [29]. As discussed below, noncontinuous decision rules are of considerable interest in economics. Finally, many applications make use of simulation-based estimation (e.g., see [14] and [22] for two recent examples); however, in general these papers do not get into the analysis of asymptotic properties of the estimator.

## 2. Assumptions.
For convenience of the presentation, the set of states $S$ is a hypercube in Euclidean space $\mathcal{R}^k$, that is, $S = \{s = (\ldots, s_i, \ldots) : a_i \leq s_i \leq b_i\}$ for given constants $a_i < b_i$ and all $i = 1, \ldots, k$ and $\geq$ is the usual Euclidean (partial) order. $\mathcal{S}$ denotes the relativized Borel $\sigma$-algebra of $S$. The vector of shocks $\varepsilon$ follows an i.i.d. process defined by a probability law $Q$ on a measurable space $(E, \mathcal{E})$. The domain of parameter vectors $\Theta \subset \mathcal{R}^l$ is a compact set. We say that a function $h$ in $S$ is increasing if $h(s) \geq h(s')$ for all $s \geq s'$.

Function $\varphi : S \times E \times \Theta \to S$ satisfies the following assumptions:

(A.1) *Measurability.* Function $\varphi : S \times E \times \Theta \to S$ is measurable on the product space $\mathcal{S} \times \mathcal{E} \times \Theta$.

(A.2) *Monotonicity.* For each $(\varepsilon, \theta)$, mapping $\varphi(\cdot, \varepsilon, \theta) : S \to S$ is increasing.

Note that no order preserving conditions are required over space $E \times \Theta$.

(A.3) *Feller property.* For every fixed $\theta$ and every continuous function $f : S \to R$,

$$(2.1) \qquad \int f(\varphi(s_j, \varepsilon, \theta))Q(d\varepsilon) \to_j \int f(\varphi(s, \varepsilon, \theta))Q(d\varepsilon) \qquad \text{as } s_j \to_j s.$$

This weak form of continuity over $S$ is usually assumed to guarantee existence of an invariant probability $\mu_\theta$ for mapping $\varphi(\cdot, \cdot, \theta)$; see [1], Chapter 1.

We also need some type of continuity on the domain of parameter vectors $\Theta$. Since mapping $\varphi(\cdot, \cdot, \theta)$ may not be continuous on $S \times E$, the continuity in $\theta$ is defined with the help of some majorizing and minorizing functions $\varphi^\kappa$ and $\varphi_\kappa$. Let $e$ denote the unit vector $(1, 1, \ldots, 1)$ and $\kappa > 0$. Then $\varphi^\kappa(s, \varepsilon, \theta) = proj_S[\varphi(proj_S[s + \kappa e], \varepsilon, \theta) + \kappa e]$ for all $(s, \varepsilon)$, where $s' = proj_S[s + \kappa e]$ means the natural projection or minimum Euclidean distance of vector $s + \kappa e$ to set $S$. The projection mapping just ensures that function $\varphi^\kappa$ is well defined. Since function $\varphi$ is monotone in $s$ and the positive vector $\kappa e$ is added to both vector $s$ in the domain and to the value $\varphi(s', \varepsilon, \theta)$, we get that $\varphi^\kappa(s, \varepsilon, \theta) \geq \varphi(s, \varepsilon, \theta)$. We write $\varphi^\kappa(\cdot, \cdot, \theta) \geq \varphi(\cdot, \cdot, \theta)$ to express functional dominance over $\geq$. In the same way, we define function $\varphi_\kappa(s, \varepsilon, \theta) = proj_S[\varphi(proj_S[s - \kappa e], \varepsilon, \theta) - \kappa e]$. Functions $\varphi^\kappa$ and $\varphi_\kappa$ will be further discussed below [Section 5, remark (iv)].

(A.4)  *Continuity in the parameters.* For every $\theta$ and $\kappa > 0$, the set $\{\theta' | \varphi^\kappa(\cdot, \cdot, \theta) \geq \varphi(\cdot, \cdot, \theta') \geq \varphi_\kappa(\cdot, \cdot, \theta)\}$ contains an open neighborhood $V_\kappa(\theta)$ of $\theta$.

Observe that (A.1)–(A.3) will all be satisfied if $\varphi$ is a continuous function over a compact domain $S \times E \times \Theta$. Moreover, (A.4) would hold trivially under the upper and lower envelope functions $\varphi^\kappa(s_0, \varepsilon, \theta) = \sup_{\theta' \in V_\kappa(\theta)} \varphi(s_0, \varepsilon, \theta')$ and $\varphi_\kappa(s_0, \varepsilon, \theta) = \inf_{\theta' \in V_\kappa(\theta)} \varphi(s_0, \varepsilon, \theta')$.

(A.5)  *Uniqueness of the invariant probability.* For every $\theta$, mapping $\varphi(\cdot, \cdot, \theta)$ has a unique invariant probability $\mu_\theta$.

Certain conditions guarantee (A.5), for example, [2] and [13]. This assumption will simplify the analysis considerably, and it is necessary to obtain our global convergence results. These results can suitably be extended to account for multiple invariant probabilities; e.g., see Lemma 4.2 below.

Let $\| \cdot \|$ be the max norm in $\mathcal{R}^l$. Then for every fixed $\theta$, we define the distance between mappings $\varphi(\cdot, \cdot, \theta)$ and $\tilde\varphi(\cdot, \cdot, \theta)$ over $S \times E$ as

$$(2.2) \quad d(\varphi(\cdot, \cdot, \theta), \tilde\varphi(\cdot, \cdot, \theta)) = \max_{s \in S} \left[ \int \| \varphi(s, \varepsilon, \theta) - \tilde\varphi(s, \varepsilon, \theta) \| Q(d\varepsilon) \right].$$

The above assumptions ensure that the invariant probability $\mu_\theta$ is always well defined and weakly continuous in $\theta$. Moreover, the maximum in (2.2) is always attained at some $s$. It should be emphasized that the standard sup norm in the space of functions $\varphi(\cdot, \cdot, \theta)$ will be very restrictive since $\varphi(\cdot, \cdot, \theta)$ may not be continuous in $(s, \varepsilon)$.

### 3. The simulated moments estimator (SME).

3.1. *Model simulation.* Let us first place ourselves in a simple framework of model simulation. Assume that a researcher can evaluate function $\varphi$ at any given point $(s, \varepsilon, \theta)$, and can draw sequences $\{\hat{\varepsilon}_n\}_{n \geq 1}$ from a random generator that mimics the distribution of $\{\varepsilon_n\}_{n \geq 1}$. No knowledge of the actual realization of the shock process $\{\varepsilon_n\}_{n \geq 1}$ is required. Later, the analysis is extended to the more typical situation in which the researcher can only perform evaluations of a numerical approximation $\varphi^j$. Hence, for each parameter value $\theta$ and initial condition $s_0$ one can generate sequences $\{s_n(s_0, \omega, \theta)\}_{n \geq 1}$ using dynamical system (1.1); that is, $s_n(s_0, \omega, \theta) = \varphi(s_{n-1}(s_0, \omega, \theta), \varepsilon_n, \theta)$ for all $n \geq 1$. It should be emphasized that in order to search for the true value $\theta^0$ we need to sample over the whole parameter space $\Theta$.

3.2. *Probability spaces.* Let $\tilde{\mathbf{s}} = \{\tilde{s}_n\}_{n \geq 1}$ be a sample path of observations of the data generating process. Let $\omega = \{\varepsilon_n\}_{n \geq 1}$ be a sequence of realizations of the shock process. A measure $\tilde{\gamma}$ is defined over the space of sequences $\tilde{\mathbf{s}} = \{\tilde{s}_n\}_{n \geq 1}$. Also, from the probability law $Q$ a measure $\gamma$ can be constructed over the space $\Omega$ of sequences $\omega = \{\varepsilon_n\}_{n \geq 1}$. Let $\lambda = \gamma \times \tilde{\gamma}$ represent the product measure.

3.3. *The SME.* Our definition of the SME is conformed by the following elements:

(1) A *function of interest* $f : S \to \mathcal{R}^p$ assumed to be continuous and increasing. Then the expected or mean value $E_\theta(f) = \int f(s)\mu_\theta(ds)$ may represent $p$ moments of an invariant distribution $\mu_\theta$. The expected value of $f$ over the invariant distribution of the data generating process will be denoted by $\bar{f}$. The restriction to increasing functions is rather harmless. Indeed, following [32] the subclass of continuous and increasing functions determines convergence in the weak topology of measures; that is, a sequence of probability measures $\{\mu_n\}_{n \geq 1}$ converges weakly to measure $\mu$ if and only if $\int f(s)\mu_n(ds) \to_n \int f(s)\mu(ds)$ for every continuous and increasing function $f : S \to \mathcal{R}$.

(2) A *distance function* $G : \mathcal{R}^p \times \mathcal{R}^p \to \mathcal{R}$ assumed to be continuous. The minimum distance is attained at a vector of parameter values

$$(3.1) \qquad \theta^0 = \arg\min_{\theta \in \Theta} G(E_\theta(f), \bar{f}).$$

Under the above assumptions, there exists an optimal solution $\theta^0$. To facilitate the presentation we suppose that $\theta^0$ is unique.

(3) An *estimation rule* characterized by a sequence of distance functions $\{G_N\}_{N \geq 1}$ and choices for the horizon $\{\tau_N\}_{N \geq 1}$ of the model's simulations. This rule yields a sequence of estimated values $\{\hat{\theta}_N\}_{N \geq 1}$ from associated optimization problems with finite samples of model's simulations and data. The estimated value $\hat{\theta}_N(s_0, \omega, \tilde{\mathbf{s}})$ is obtained as

$$
\hat{\theta}_N(s_0, \omega, \tilde{\mathbf{s}}) = \arg\inf_{\theta \in \Theta} G_N\left(\frac{1}{\tau_N(\omega, \tilde{\mathbf{s}})} \sum_{n=1}^{\tau_N(\omega, \tilde{\mathbf{s}})} f(s_n(s_0, \omega, \theta)),\right.
$$
(3.2)
$$
\left.\frac{1}{N}\sum_{n=1}^{N} f(\tilde{s}_n), \omega, \tilde{\mathbf{s}}\right).
$$

We assume that the sequence of continuous functions $\{G_N(\cdot, \cdot, \omega, \tilde{\mathbf{s}})\}_{N \geq 1}$ converges uniformly to function $G(\cdot, \cdot)$ for $\lambda$-almost all $(\omega, \tilde{\mathbf{s}})$, and the sequence of functions $\{\tau_N(\omega, \tilde{\mathbf{s}})\}_{N \geq 1}$ goes to $\infty$ for $\lambda$-almost all $(\omega, \tilde{\mathbf{s}})$. Note that both functions $G_N(\cdot, \cdot, \omega, \tilde{\mathbf{s}})$ and $\tau_N(\omega, \tilde{\mathbf{s}})$ are allowed to depend on $\omega$ and $\tilde{\mathbf{s}}$. These functions will usually depend on all information available up to time $N$. The rule $\tau_N$ reflects that model's simulations may be of a different length than data samples.

It should be stressed that problem (3.1) is defined over population characteristics of the model and of the data generating process, whereas problem (3.2) is defined over statistics of finite simulations and data.

DEFINITION. The SME is a sequence of measurable functions $\{\hat{\theta}_N(s_0, \omega, \tilde{\mathbf{s}})\}_{N \geq 1}$ such that each function $\hat{\theta}_N$ satisfies (3.2) for all $s_0$ and $\lambda$-almost all $(\omega, \tilde{\mathbf{s}})$.

By the measurable selection theorem ([8], Chapter 2), there exists a sequence of measurable functions $\{\hat{\theta}_N\}_{N \geq 1}$.

## 4. Main results.
As stressed above, the proof of consistency of the SME has to deal with a continuum of invariant probabilities. The strategy is to show that minimization problem (3.1) can be approximated by a sequence of optimization problems (3.2) for $\lambda$-almost all $(\omega, \tilde{\mathbf{s}})$. The hardest step in the proof is to demonstrate the uniform convergence (a.s.) of the sequences $\{\frac{1}{\tau_N(\omega,\tilde{\mathbf{s}})}\sum_{n=1}^{\tau_N(\omega,\tilde{\mathbf{s}})} f(s_n(s_0,\omega,\theta))\}_{N \geq 1}$ to $E_\theta(f)$ over $\Theta$, as $N \to \infty$. This uniform convergence of the simulated statistics follows from a *sandwich* argument that builds upon the weak continuity of the invariant probabilities of functions $\varphi^\kappa$ and $\varphi_\kappa$ under perturbations in $\kappa$ (Lemma 4.1), a generalized law of large numbers for each individual function (Lemma 4.2), and the order preserving property of $\varphi$ as stated in (A.2). These lemmas draw on some results in [30] on *pointwise* convergence properties of invariant probabilities for random dynamical systems satisfying the Feller property. We extend these results to establish *uniform* convergence over $\Theta$ of the sample distributions (Proposition 4.3). Proofs are gathered in the final section.

Let $\mu_\theta^\kappa$ be an invariant probability of function $\varphi^\kappa(\cdot, \cdot, \theta)$, and $\mu_{\kappa\theta}$ an invariant probability of function $\varphi_\kappa(\cdot, \cdot, \theta)$.

LEMMA 4.1 (Continuity of the invariant probabilities). *Assume that the sequence of positive scalars $\{\kappa_j\}_{j\geq 1}$ converges to zero. Then under* (A.1)–(A.5) *every sequence of invariant probabilities $\{\mu_\theta^{\kappa_j}\}_{j\geq 1}$ (resp. $\{\mu_{\kappa_j\theta}\}_{j\geq 1}$) converges weakly to the unique invariant probability $\mu_\theta$ of mapping $\varphi(\cdot, \cdot, \theta)$.*

Note that our primitive assumptions (A.1)–(A.5) do not rule out the possibility that the auxiliary functions $\varphi^\kappa(\cdot, \cdot, \theta)$ and $\varphi_\kappa(\cdot, \cdot, \theta)$ may contain multiple invariant probabilities. Hence, let $\mathcal{I}(\mu_\theta^\kappa)$ be the set of all the invariant probabilities $\mu_\theta^\kappa$ of $\varphi^\kappa(\cdot, \cdot, \theta)$. From the proof of Lemma 4.1, the set $\mathcal{I}(\mu_\theta^\kappa)$ is compact and convex in the weak topology of measures. Then every continuous linear functional $\mu_\theta^\kappa \to \int f(s)\mu_\theta^\kappa(ds)$ attains a maximum and a minimum over the domain $\mathcal{I}(\mu_\theta^\kappa)$. In the same way, let $\mathcal{I}(\mu_{\kappa\theta})$ be the set of all the invariant probabilities $\mu_{\kappa\theta}$ of $\varphi_\kappa(s, \varepsilon, \theta)$. Finally, for every sequence of shocks $\omega = \{\varepsilon_n\}_{n\geq 1}$ and initial condition $s_0$, define recursively the sample paths $\{s_n^\kappa(s_0, \omega, \theta)\}_{n\geq 1}$ and $\{s_{\kappa n}(s_0, \omega, \theta)\}_{n\geq 1}$ generated by functions $\varphi^\kappa$ and $\varphi_\kappa$ as $s_n^\kappa(s_0, \omega, \theta) = \varphi^\kappa(s_{n-1}^\kappa(s_0, \omega, \theta), \varepsilon_n, \theta)$ and $s_{\kappa n}(s_0, \omega, \theta) = \varphi_\kappa(s_{\kappa n-1}(s_0, \omega, \theta), \varepsilon_n, \theta)$, for all $n \geq 1$.

We next show that the range of variation of the average behavior of a typical simulated path $\{s_n^\kappa(s_0, \omega, \theta)\}_{n\geq 1}$ or $\{s_{\kappa n}(s_0, \omega, \theta)\}_{n\geq 1}$ is bounded by the corresponding expected values over the sets of invariant probabilities $\mathcal{I}(\mu_\theta^\kappa)$ and $\mathcal{I}(\mu_{\kappa\theta})$.

LEMMA 4.2 (A generalized law of large numbers). *Under* (A.1)–(A.5), *for every fixed $\theta$ in $\Theta$,*

$$
\begin{aligned}
\max_{\mu_\theta^\kappa \in \mathcal{I}(\mu_\theta^\kappa)} \int f(s)\mu_\theta^\kappa(ds) &\geq \limsup_{N\to\infty} \frac{1}{N}\sum_{n=1}^N f(s_n^\kappa(s_0, \omega, \theta)) \\
&\geq \liminf_{N\to\infty} \frac{1}{N}\sum_{n=1}^N f(s_{\kappa n}(s_0, \omega, \theta)) \\
&\geq \min_{\mu_{\kappa\theta} \in \mathcal{I}(\mu_{\kappa\theta})} \int f(s)\mu_{\kappa\theta}(ds)
\end{aligned}
\tag{4.1}
$$

*for all $s_0$ and $\gamma$-almost all $\omega$.*

In other words, there exists a set $\hat{\Omega}$ with $\gamma(\hat{\Omega}) = 1$ such that (4.1) is satisfied for all $(s_0, \omega) \in S \times \hat{\Omega}$. If each of the functions $\varphi^\kappa(\cdot, \cdot, \theta)$ and $\varphi_\kappa(\cdot, \cdot, \theta)$ has a unique invariant distribution, then Lemma 4.2 is a simple consequence of the law of large numbers of [4] together with (A.2).

The foregoing lemmas are the main ingredients in the proof of the following uniform law of large numbers over the parameter space $\Theta$. This result is key to substantiate the strong consistency of the SME.

PROPOSITION 4.3 (Uniform convergence of the simulated statistics). *Under* (A.1)–(A.5), *there is a set* $\hat{\Omega}$ *with* $\gamma(\hat{\Omega}) = 1$ *such that every pair* $(s_0, \omega) \in S \times \hat{\Omega}$ *satisfies the following property: For each* $\epsilon > 0$, *there is a constant* $N_\epsilon(\omega)$ *such that*

$$(4.2) \qquad \left| \frac{1}{N} \sum_{n=1}^{N} f(s_n(s_0, \omega, \theta)) - E_\theta(f) \right| < \epsilon$$

*for all* $N \geq N_\epsilon(\omega)$ *and all* $\theta$ *in* $\Theta$.

Note that $N_\epsilon(\omega)$ only depends on $\omega$, and hence it holds for all $s_0$. This proposition could be restated in terms of the uniform convergence in $\theta$ of the sample distributions (e.g., see [20]), as the set of continuous functions $f : S \to \mathcal{R}^p$ is separable.

THEOREM 4.4 (Strong consistency of the SME). *Assume that the process* $\tilde{\mathbf{s}} = \{\tilde{s}_n\}_{n \geq 0}$ *is stationary and ergodic. Then under* (A.1)–(A.5), *for all* $s_0$ *and* $\lambda$-*almost all* $(\omega, \tilde{\mathbf{s}})$, *the SME* $\{\hat{\theta}_N(s_0, \omega, \tilde{\mathbf{s}})\}_{N \geq 1}$ *converges to* $\theta^0$.

The SME is computationally costly and extensive model simulations must be performed to sample the region of feasible parameter values. A gain in computational efficiency should be obtained when some parameter values are constrained to take on certain values or can be estimated by more practical methods. For instance, let $\theta = (\theta_1, \theta_2)$ and suppose that the second component $\theta_2$ can be estimated by traditional methods. Then similar consistency results can be established for vector $\theta_1$ under the presumption that the estimator for vector $\theta_2$ is consistent.

The consistency of the estimator can also be extended to numerical approximations. In most dynamical models, the solution $\varphi$ does not admit an analytical representation, but it is usually possible to perform functional evaluations of a numerical approximation. Most numerical methods can generate sequences of solutions $\{\varphi^j\}_{j \geq 1}$ that converge to the original function $\varphi$ as we refine the approximation scheme. Hence, it is of interest to know asymptotic properties of the estimator as the numerical error vanishes.

As in (3.1), for every approximate function $\varphi^j$ a solution $\theta^j$ is attained over the parameterized family of invariant probabilities. More specifically,

$$(4.3) \qquad \theta^j = \arg\min_{\theta \in \Theta} G\left( \int f(s) \mu_\theta^j(ds), \bar{f} \right),$$

where $\mu_\theta^j$ is an invariant probability of mapping $\varphi^j(\cdot, \cdot, \theta)$. The invariant probability of $\mu_\theta^j$ may not be unique, even though for each $\theta$ the original mapping $\varphi(\cdot, \cdot, \theta)$ is assumed to have a unique invariant probability $\mu_\theta$. Also, the solution $\theta^j$ may not be unique. The idea is that certain primitive assumptions (cf. [2] and [13]) may guarantee uniqueness of the invariant distribution $\mu_\theta$ of $\varphi(\cdot, \cdot, \theta)$ but uniqueness is not generally preserved under numerical perturbations of the model such as piecewise linear or polynomial interpolations. Hence, problem (4.3) should be understood as a minimization over the correspondence of invariant distributions $\mu_\theta^j$.

THEOREM 4.5 (Consistency of the SME for numerical approximations). *Let $\varphi^j$ satisfy (A.1) and (A.3) for every $j$. Let $\varphi$ satisfy (A.1) and (A.3)–(A.5). Assume that the sequence of functions $\{\varphi^j(\cdot, \cdot, \theta)\}_{j \geq 1}$ converges uniformly to $\varphi(\cdot, \cdot, \theta)$ over $\Theta$ in the metric (2.2). Then every sequence of optimal solutions $\{\theta^j\}_{j \geq 1}$ of (4.3) must converge to the original solution $\theta^0$ of (3.1).*

Obviously, Theorems 4.4 and 4.5 can be combined to obtain convergence of the estimator in both $N$ and $j$. Note that Theorem 4.5 does not depend on monotonicity condition (A.2) since we are assuming the uniform convergence of the sequence of functions $\{\varphi^j(\cdot, \cdot, \theta)\}_{j \geq 1}$. But as in the proof of Proposition 4.3, condition (A.2) can be instrumental to build a sandwich argument to establish the uniform convergence of $\{\varphi^j(\cdot, \cdot, \theta)\}_{j \geq 1}$ over $\theta$ in $\Theta$ from the pointwise convergence of these functions for each $\theta$.

## 5. Remarks.
The following issues may deserve further discussion:

(i) *Latent variables.* The state space $S$ may contain both observable and unobservable state variables. For instance, assume that $s = (s_1, s_2)$ is a vector in $\mathcal{R}^2$ such that $s_1$ represents the value of production or output, $s_2$ represents the level of efficiency or productivity, and $\varepsilon$ is a productivity shock. Usually, both $s_2$ and $\varepsilon$ are unobservable. But function $f$ may encompass enough moments or characteristics of variable $s_1$ so as to identify the whole vector of parameters $\theta$. Several papers (e.g., see [24, 28, 33] for some early examples) have stressed the importance of simulation-based estimation for models with unobservable or so-called "latent variables."

(ii) *Monotonicity.* Our results on the consistency of the SME follow from relatively simple assumptions that are easy to check in applications. The monotonicity condition (A.2) plays a key role in our arguments, and it is the most restrictive assumption.[4] All other regularity assumptions are dictated by technical considerations and are less limiting in applications. As compared to [11], we dispense with some strong assumptions such as geometric ergodicity, a global modulus of continuity condition, and damping conditions on the dynamics of the system. All these assumptions may be hard to check in applications.

---

[4]Monotonicity can be weakened under specific functional forms. Let $s = (s_1, s_2)$ and $\theta = (\theta_1, \theta_2)$. Assume that $s_1 = \varphi_1(s_1, s_2, \theta)$ and $s_2 = \varphi_2(s_2, \theta_2)$. Now, if there is an unbiased

(iii) *Ergodic processes.* Our results could be extended beyond Markov processes, but stronger continuity conditions on mapping $\varphi$ may be required. For instance, [8], Proposition 6.21, derives a generalized law of large numbers that it is suitable for extensions of Lemma 4.2 to ergodic processes under the more restrictive assumption of continuity of mapping $\varphi(\cdot, \varepsilon, \theta) \colon S \to S$.

(iv) *Continuity in the vector of parameters.* Although function $\varphi(s, \varepsilon, \theta)$ may not be continuous in $s$, our method of proof relies on some type of continuity of $\varphi(s, \varepsilon, \theta)$ in $\theta$. A natural approach would be to require the continuity of $\varphi(s, \varepsilon, \theta)$ in $\theta$ under the distance function $d(\varphi(\cdot, \cdot, \theta), \varphi(\cdot, \cdot, \theta'))$ in (2.2). This metric is suitable for our purposes, since by (A.3) discontinuities of $\varphi(s, \varepsilon, \theta)$ in $s$ can be smoothed out when integrating over $\varepsilon$. A related notion of continuity is assumed in (A.4) under the majorizing and minorizing functions $\varphi^{\kappa}$ and $\varphi_{\kappa}$. This construction has been useful to circumvent some measurability problems. To motivate the definition of functions $\varphi^{\kappa}$ and $\varphi_{\kappa}$ the following simple example may be helpful. Suppose that $s$, $\varepsilon$ and $\theta$ are scalars. Let $\varphi(s, \varepsilon, \theta) = s + \varepsilon + \theta$ if $s + \varepsilon + \theta \leq 2$, and $\varphi(s, \varepsilon, \theta) = s + \varepsilon + \theta + 5$ if $s + \varepsilon + \theta > 2$. That is, this function is increasing, and has a jump after reaching a certain threshold. (Observe that a sufficient condition for the Feller property to be satisfied is that $\varepsilon$ has a continuous density.) Consider now the majorizing function $\hat{\varphi}^{\kappa}(s, \varepsilon, \theta) = \varphi(s, \varepsilon, \theta) + \kappa e$. Then $\hat{\varphi}^{\kappa} \geq \varphi$. However, the set $\{\theta' | \hat{\varphi}^{\kappa}(\cdot, \cdot, \theta) \geq \varphi(\cdot, \cdot, \theta') \geq \varphi(\cdot, \cdot, \theta)\}$ is empty. Note that the majorizing function in (A.4) is $\varphi^{\kappa}(s, \varepsilon, \theta) = \varphi(s', \varepsilon, \theta) + \kappa e$ where $s' = s + \kappa e$. Here, the set $\{\theta' | \varphi^{\kappa}(\cdot, \cdot, \theta) \geq \varphi(\cdot, \cdot, \theta') \geq \varphi(\cdot, \cdot, \theta)\}$ contains an open set of parameters $\theta'$.

(v) *Feller property.* [2] and [10] dispense with the Feller property by introducing a mild mixing condition that guarantees existence, uniqueness and global stability to the invariant probability. Under the mixing condition the random dynamical system is a contraction mapping in a suitable metric space of probabilities. The resulting metric topology, however, is too fine to allow for continuous perturbations of the parameter space. For instance, when applied to a probability measure $\mu$ the above restricted perturbation $\hat{\varphi}^{\kappa}(s, \varepsilon, \theta) = \varphi(s, \varepsilon, \theta) + \kappa e$ may not vary continuously with $\kappa$ over the distance function (2.4) of [2]. Therefore, in the absence of further specific conditions (e.g., [11]) the Feller property seems indispensable for the strong consistency of our simulation-based estimator.

---

estimator for $\theta_2$, consistency of the SME for $\theta_1$ can be insured by monotonicity of $\varphi_1$ in $s_1$. In the neoclassical growth model, $s_2 = \varphi_2(s_2, \theta_2)$ corresponds to the law of motion of the exogenous technological progress, and monotonicity of $s_1 = \varphi_1(s_1, s_2, \theta)$ in $s_1$ follows from the concavity of the utility and production functions.

**6. Applications.** Several dynamic optimization problems generate noncontinuous, monotone decision rules. Simple discontinuities for monotone Markov processes are often encountered in models of economic growth and finance, models of firm entry, patent races, replacement of durable goods, job search, marriage and fertility; e.g., see [3, 17, 23], and the aforementioned papers in (i) of the previous section. In these models, the optimal decision can feature isolated jumps because of discrete choices or lack of convexity but the Feller property may nevertheless be satisfied.

Our purpose here is to estimate some critical parameters of a simplified version of the stock market model of [18]. This model borrows several elements from [7] and [26] who are concerned with the effects of technology adoption on economic growth and business fluctuations rather than on financial markets. None of these three papers consider model estimation, and simply report simulations for various benchmark calibrations of the model.

6.1. *The model.* The representative household has preferences over consumption $c$ of an aggregate good and desutility of work $l$ as represented by the expected discounted objective:

$$
(6.1) \qquad E_0 \left\{ \sum_{t=0}^{\infty} \beta^t \left[ \ln(c_t) - \frac{l_t^{1+\chi}}{1+\chi} \right] \right\}
$$

with $0 < \beta < 1$ and $\chi > 0$. For given initial quantity of the aggregate asset, $\hat{a}$, the optimization problem faced by this agent is to choose a stochastic sequence of consumption, labor and shares of the aggregate stock $\{c_t, l_t, a_t\}_{t \geq 0}$ that maximizes the objective in (6.1) subject to the sequence of budget constraints

$$
(6.2) \qquad c_t + q_t a_t = w_t l_t + (q_t + d_t) a_{t-1}
$$

with $a_t \geq 0$. Observe that $q_t$ denotes the market price of one stock unit, $d_t$ denotes the dividend, and $w_t$ is the wage, for $t = 0, 1, \ldots$.

There is a mass of $A_t$ intermediate goods that are bundled together into a composite good $M_t$ defined by a CES technology, $M_t = [\int_0^{A_t} m_t(j)^{1/\vartheta} \, dj]^\vartheta$ where $m(j)$ denotes the amount of intermediate good $j \in [0, A_t]$ bought by the firm and $\vartheta > 1$. The firm producing the final good accumulates capital $k$ and buys labor $l$ and a bundle of intermediate goods $M$ using production function $y = \theta_t (k_t^\alpha l_t^{1-\alpha})^{1-\gamma} M_t^\gamma$. Output can be devoted to consumption or investment, and capital depreciates at a rate $\delta > 0$. Total factor productivity $\theta_t$ is governed by the following law of motion: $\log \theta_t = \varphi_\theta \ln \theta_{t-1} + \varepsilon_t^\theta$, for $\varepsilon_t^\theta \sim N(0, \sigma_\theta)$. The firm chooses stochastic sequences of investment, labor and intermediate goods $\{i_t, l_t, m_t(j)_{j \in [0, A_t]}\}_{t \geq 0}$ to maximize expected discounted revenues.

Monopolistic competition prevails in the intermediate goods market. That is, producer of intermediate good $j$ selects both the optimal quantity $m(j)$ to sell and corresponding price $p(j)$—taken as given prices and quantities set up by all other producers of intermediate goods. Following [26], we postulate a very simple technology: Production of one unit of good $j$ just requires one unit of the final good. Then, the profit at time $t$ of the producer of variety $j$ is given by

$$(6.3) \qquad \pi_t(j) \equiv \max_{m_t(j)} \{ p_t(j) m_t^i(j) - m_t^i(j) \}.$$

Without loss of generality, we suppose that $\pi_t(j)$ is the same for all $j$. Each intermediate good may eventually be displaced from the market. Let $\phi$ be the probability of survival of good $j$. Let $r_t$ be the stochastic rate of interest at time $t$. Then for $\eta_t^s = (1 + r_{t+1}) \cdots (1 + r_s)$ the present value $V_t$ of operating each technology from the beginning of time $t$ is defined as:

$$(6.4) \qquad V_t = E_t \left\{ \sum_{s=t+1}^{\infty} \frac{\pi_s(j)}{\eta_t^s} \phi^{s-t} \right\}.$$

Technological innovations arrive exogenously to the economy. The total stock of technological innovations $Z_t$ evolves according to the law of motion

$$(6.5) \qquad Z_t = \phi Z_{t-1} + x_t$$

with

$$(6.6) \qquad \ln x_t = \varphi_x \ln x_{t-1} + \varepsilon_t^x, \qquad \varphi_x \in (0,1), \qquad \varepsilon_t^x \sim N(0, \sigma_x).$$

These functional forms are considered here for simplicity. Indeed, we could allow for discontinuities in variables $Z$ and $x$ to reflect sudden changes in the transmission of new technologies. These technologies are put into use by local adopters. The adoption sector behaves as a perfectly competitive market. Let $A_t$ be the stock of already adopted technologies, and $\lambda(H)$ the probability of adoption of a new technology after investing the amount of resources $H$. Then

$$(6.7) \qquad A_{t+1} = \lambda(H_t)[Z_t - A_t] + \phi A_t.$$

The optimal amount of expenditure $H_t$ is derived from the following Bellman equation in which the value function is the option value $J_t$ of a new technology.[5]

$$(6.8) \quad J_t = \max_{H_t} \left\{ -H_t + \phi E_t \left[ \frac{1}{1 + r_{t+1}} (\lambda(H_t) V_{t+1} + (1 - \lambda(H_t)) J_{t+1}) \right] \right\}.$$

---

[5]As is well known, this equation can be computed recursively by the method of successive approximations.

6.2. *Equilibrium and asset prices.* In our model, the exogenous state variables are the stock of technological innovations $Z_t$ and $x_t$ and the value of total factor productivity $\theta_t$, and the endogenous state variables are the amount of capital $k_t$ and the stock of adopted technologies $A_t$. The remaining variables are determined endogenously as solutions to the model under (6.1)–(6.8) and the market clearing conditions. Let us assume that $a_t = 1$ so that $q_t$ corresponds to the value of the stock market. This value can be decomposed into the value of existing stocks plus the option value of current and future technological innovations.

PROPOSITION 6.1. *The stock market value $q_t$ can be computed as*

$$(6.9) \qquad q_t = k_{t+1} + V_t A_t + J_t^+(Z_t - A_t) + \xi_t,$$

*where $J_t^+ \equiv J_t + H_t$, and $\xi_t \equiv E_t\{\sum_{s=t+1}^{\infty} \frac{1}{\eta_t^s} J_s(Z_s - \phi Z_{s-1})\}$.*

Hence, the value of the stock market is given by the sum of four components: The replacement cost of installed capital, the value of adopted technologies, the option value of technological innovations currently available but not yet implemented, and the present value of technological innovations expected to occur. Most economic models identify the stock market value $q_t$ with the capital stock $k_{t+1}$.

6.3. *Computation, calibration and estimation of the model.* The model can be solved numerically from the first-order conditions of the above optimization problems and the market clearing conditions. We compute these equations using a low degree perturbation method; see [18]. Various parameters of the model are calibrated to match some empirical statistics of medium-term fluctuations and the volatility of patents.[6] But for reasons that will become clear below, we apply our simulation-based estimator to three critical parameters: $\chi$, $\gamma$ and $\vartheta$.

Following the economics literature (e.g., see [7] and references therein) we choose values for the set of parameters $(\beta, \alpha, \delta)$ to match US macro data.[7] Parameter $\beta$ is fixed at 0.95 so that the real rate of return of investment in the deterministic steady state is equal to the average real return on the S&P index over the period 1948–2004. Parameter $\alpha$ is set at 0.3 to match the average income share of labor costs, and the annual depreciation rate $\delta$ is set at 0.075.

---

[6] Following [7] medium-term cycles are defined as those within a frequency band of 2–50 years. The data are filtered using the band-pass filter of [5].

[7] We consider output, hours, labor productivity and TFP for the nonfarm business sector. The source is the Bureau of Labor Statistics (BLS). Consumption is measured as the sum of nondurables and services and investment is nonresidential. Both series are obtained from the Bureau of Economic Analysis (BEA). Each variable is transformed in per capita terms using the population aged 15–64. The data sample spans from 1948 to 2004.

As a proxy for the number of adopted technologies $A_t$ we use the percentage deviation from its exponential trend of the number of patents issued.[8] The parameter values for the law of motion of variables $\theta$ and $Z$ are selected to approximate the variance, correlation, and first-order autocorrelation of total factor productivity and patents over medium-term cycles in the data.

Following [7], we set the value of $\phi$ to 0.97 and the steady state value of $\lambda$ to 0.1. We assume that the probability of adoption takes on the functional form

$$(6.10) \qquad \lambda_t(H_t) = \Lambda \left( \frac{A_t}{k_t} H_t \right)^\rho$$

with constant $\Lambda > 0$ and $\rho \in (0, 1)$. We estimate that the average share of adoption expenditures over sales[9] is 0.019 for the period 1970–1998. Then, in the steady state solution, $\rho$ must be 0.11.

We are then left with the estimation of parameters $\chi$, $\gamma$ and $\vartheta$. There are several estimates of the elasticity of individual labor supply $\chi$, but these estimates do not seem adequate for our aggregate model since hours worked fluctuate much less than in the data. Regarding the share of intermediate goods $\gamma$, there are empirical estimates for the industrial sector (e.g., [19]), but our measure of intermediate goods is much broader. Along the same lines, there are estimates of the mark-up parameter $\vartheta$ for various sectors [27], but our model includes a very stylized cost function for the production of intermediate goods and there are various nontangible intermediate goods that may be expensed as patents and trademarks.

We define the objective function of our estimator as

$$(6.11) \qquad \frac{1}{\hat{\sigma}_{\hat{\sigma}_{inv}}} (\sigma_{inv} - \hat{\sigma}_{inv})^2 + \frac{1}{\hat{\sigma}_{\hat{\sigma}_{hours}}} (\sigma_{hours} - \hat{\sigma}_{hours})^2$$
$$+ \frac{1}{\hat{\sigma}_{\hat{\sigma}_{stock}}} (\sigma_{stock} - \hat{\sigma}_{stock})^2,$$

where the standard deviations of the model $\sigma$ are compared against standard deviations of the data $\hat{\sigma}$ for variables $inv$ = investment, $hours$ = hours worked, and $stock$ = stock market value. These deviations are weighted by the inverse of the standard deviation of $\hat{\sigma}$ of each corresponding variable. Note that each $\sigma$ is a function of parameter values. The objective is minimized over the Euclidean product of interval values: $\chi \in [0, 1]$, $\gamma \in [0.1, 0.60]$ and $\vartheta \in [1.1, 2.2]$ under $\hat{\sigma}_{inv} = 8.86$, $\hat{\sigma}_{hours} = 3.31$, $\hat{\sigma}_{stock} = 31.41$, and $\hat{\sigma}_{\hat{\sigma}_{inv}} = 0.0091$, $\hat{\sigma}_{\hat{\sigma}_{hours}} = 0.0035$, and $\hat{\sigma}_{\hat{\sigma}_{stock}} = 0.0315$.

---

[8]The data come from the total number of utility patents granted from the US Patent and Trademark Office for 1970–2004, and from Historical Statistics of the United States series W-99 for 1948–1970.

[9]This measure is estimated as the average ratio of nonfederally funded development expenditures (a subset of R&D expenditures) over net sales for R&D-performing companies from the National Science Foundation.

The minimum is achieved at the following parameters: $\chi = 0$, $\gamma = 0.322$ and $\vartheta = 1.4$, with values for the objective: $\sigma_{inv} = 10.44$, $\sigma_{hours} = 1.46$ and $\sigma_{stock} = 5.34$. Therefore, one main conclusion from this exercise is that in this model the arrival of new technologies can at most account for approximately one sixth of the variability in the stock market (i.e., $\sigma_{stock} = 5.34$ vs. $\hat{\sigma}_{stock} = 31.41$), whereas we have checked that standard versions of the neoclassical model can only account for approximately one tenth of this variability. Also, the model cannot account for over half of the volatility of hours ($\sigma_{hours} = 1.46$ vs. $\hat{\sigma}_{hours} = 3.31$), which of course may suggest that some labor market frictions must play an important role.

To have a better view of the variability of stock markets, in future research we are planning to consider some other variables such as debt (leverage), liquidity constraints, taxes and subsidies, other production functions and markup policies, and monetary and financial shocks. Of course, the volatility of stock markets at present times may be due to lack of confidence and collateral requirements, but these latter variables may have played a much smaller role in our period of analysis.

## 7. Proofs.

PROOF OF LEMMA 4.1. Let $\tilde{\varphi}^{\kappa_j}(s, \varepsilon, \theta) = \varphi(s', \varepsilon, \theta)$ for $s' = proj_S[s + \kappa_j e]$ for every $\kappa_j > 0$ for $j = 1, 2, \ldots$. Let $\Psi^{\kappa_j}(s, \theta) = \int \tilde{\varphi}^{\kappa_j}(s, \varepsilon, \theta) Q(d\varepsilon)$ and $\Psi(s, \theta) = \int \varphi(s, \varepsilon, \theta) Q(d\varepsilon)$. Then

$$(7.1) \qquad \Psi^{\kappa_j}(s, \theta) - \Psi(s, \theta) = \int [\tilde{\varphi}^{\kappa_j}(s, \varepsilon, \theta) - \varphi(s, \varepsilon, \theta)] Q(d\varepsilon).$$

By (A.3), mapping $\Psi(\cdot, \theta)$ is continuous in $s$. As $S$ is a compact set, $\Psi(\cdot, \theta)$ is uniformly continuous. Hence, the sequence of functions $\{\Psi^{\kappa_j}(\cdot, \theta)\}_{j \geq 1}$ must converge uniformly to $\Psi(\cdot, \theta)$ over $S$ as $\kappa_j$ goes to zero. Further, by (A.2), we have that $\tilde{\varphi}^{\kappa_j}(s, \varepsilon, \theta) \geq \varphi(s, \varepsilon, \theta)$ for all $s$ and $\kappa_j > 0$. Then from (7.1) and the aforementioned convergence of the sequence $\{\Psi^{\kappa_j}(\cdot, \theta)\}_{j \geq 1}$ the sequence $\{\tilde{\varphi}^{\kappa_j}(\cdot, \cdot, \theta)\}_{j \geq 1}$ must converge to $\varphi(\cdot, \cdot, \theta)$ in the metric (2.2) as $\kappa_j$ goes to zero. Since $\varphi^{\kappa_j}(s, \varepsilon, \theta) = \tilde{\varphi}^{\kappa_j}(s, \varepsilon, \theta) + \kappa_j e$, the sequence of functions $\{\varphi^{\kappa_j}(\cdot, \cdot, \theta)\}_{j \geq 1}$ must converge to $\varphi(\cdot, \cdot, \theta)$ in the metric (2.2). Therefore, the corresponding sequence of invariant probabilities $\{\mu_\theta^{\kappa_j}\}_{j \geq 1}$ must converge weakly to $\mu_\theta$; see [30], Theorem 2. By the same argument, we can establish the weak convergence of the sequence of invariant probabilities $\{\mu_{\kappa_j\theta}\}_{j \geq 1}$. The proof is complete. $\square$

PROOF OF LEMMA 4.2. From the proof of the preceding lemma, it becomes clear that both functions $\varphi^\kappa$ and $\varphi_\kappa$ satisfy (A.3). Then the first and third inequalities in (4.2) follow from [30], Theorem 3. The second inequality follows from (A.2) and the fact that $\varphi^\kappa \geq \varphi_\kappa$. □

PROOF OF PROPOSITION 4.3. Since a countable union of sets of measure zero has measure zero, it suffices to establish that for a fixed rational number $\epsilon > 0$ there is a set $\hat{\Omega}$ with $\gamma(\hat{\Omega}) = 1$ such that for every $s_0$ and $\omega \in \hat{\Omega}$ we can find $N_\epsilon(\omega)$ so that

$$(7.2) \qquad \left| \frac{1}{N} \sum_{n=1}^{N} f(s_n(s_0, \omega, \theta)) - E_\theta(f) \right| < \epsilon$$

for all $N \geq N_\epsilon(\omega)$ and all $\theta$ in $\Theta$.

As $\Theta$ is compact, by Lemma 4.1 we can cover this set by a finite number of open neighborhoods $V_{\kappa_j}(\theta_j)$ for points $\{\theta_j\}$ such that for all $j = 1, \ldots, J$,

$$(7.3) \qquad \max_{\mu_{\theta_j}^{\kappa_j} \in \mathcal{I}(\mu_{\theta_j}^{\kappa_j})} \int f(s) \mu_{\theta_j}^{\kappa_j}(ds) - \min_{\mu_{\kappa_j \theta_j} \in \mathcal{I}(\mu_{\kappa_j \theta_j})} \int f(s) \mu_{\kappa_j \theta_j}(ds) < \frac{\epsilon}{2}.$$

By (A.2) and the definition of the functions $\varphi^{\kappa_j}(s, \varepsilon, \theta_j)$ and $\varphi_{\kappa_j}(s, \varepsilon, \theta_j)$, for all $\theta \in V_{\kappa_j}(\theta_j)$ and all $N \geq 1$ the following inequalities must hold true

$$(7.4) \qquad \begin{aligned} \frac{1}{N} \sum_{n=1}^{N} f(s_n^{\kappa_j}(s_0, \omega, \theta_j)) &\geq \frac{1}{N} \sum_{n=1}^{N} f(s_n(s_0, \omega, \theta)) \\ &\geq \frac{1}{N} \sum_{n=1}^{N} f(s_{\kappa_j n}(s_0, \omega, \theta_j)) \end{aligned}$$

and

$$(7.5) \qquad \int f(s) \mu_{\theta_j}^{\kappa_j}(ds) \geq \int f(s) \mu_\theta(ds) \geq \int f(s) \mu_{\kappa_j \theta_j}(ds)$$

for all invariant distributions $\mu_{\theta_i}^{\kappa_j}$ and $\mu_{\kappa_j \theta_j}$. Moreover, by Lemma 4.2 there is a set $\hat{\Omega}_j$ with $\gamma(\hat{\Omega}_j) = 1$ such that for each $(s_0, \omega) \in S \times \hat{\Omega}_j$ and $\frac{\epsilon}{2}$ there is $N_\epsilon^j(\omega)$ such that

$$(7.6) \qquad \max_{\mu_{\theta_j}^{\kappa_j} \in \mathcal{I}(\mu_{\theta_j}^{\kappa_j})} \int f(s) \mu_{\theta_j}^{\kappa_j}(ds) + \frac{\epsilon}{2} \geq \frac{1}{N} \sum_{n=1}^{N} f(s_n^{\kappa_j}(s_0, \omega, \theta_j))$$

and

$$(7.7) \qquad \frac{1}{N} \sum_{n=1}^{N} f(s_{\kappa_j n}(s_0, \omega, \theta_i)) \geq \min_{\mu_{\kappa_j \theta_j} \in \mathcal{I}(\mu_{\kappa_j \theta_j})} \int f(s) \mu_{\kappa_j \theta_j}(ds) - \frac{\epsilon}{2}$$

for all $N \geq N_\epsilon^j(\omega)$. Let $\hat{\Omega} = \bigcap_{j=1}^J \hat{\Omega}_j$ and $N_\epsilon(\omega) = \max\{N_\epsilon^j(\omega)\}_{j=1}^J$. Then, combining inequalities (7.3)–(7.7) we get that for every $s_0$ and every $\omega \in \hat{\Omega}$ condition (7.2) must hold true for all $N \geq N_\epsilon(\omega)$. □

PROOF OF THEOREM 4.4. By assumption, the process $\tilde{\mathbf{s}} = \{\tilde{s}_n\}_{n \geq 0}$ is stationary and ergodic. Hence, $\{\frac{1}{N}\sum_{n=1}^N f(\tilde{s})\}_{N \geq 1}$ converges (a.s.) to $\bar{f}$. Then Theorem 4.4 is an immediate consequence of Proposition 4.3 and the following assumptions: (i) The space of parameters $\Theta$ is a compact set, (ii) the maximizer $\theta^0$ in (3.1) is unique and (iii) the sequence of continuous functions $\{G_N(\cdot, \cdot, s_0, \omega)\}_{N \geq 1}$ converges uniformly to continuous function $G(\cdot, \cdot)$ for all $s_0$ and almost all $\omega$. □

PROOF OF THEOREM 4.5. Observe that by assumption the sequence $\{\varphi^j(\cdot, \cdot, \theta)\}_{j \geq 1}$ converges uniformly to $\varphi(\cdot, \cdot, \theta)$ in the metric $d(\cdot, \cdot)$ of (2.2). Hence, by (A.4) and a simple application of the triangle inequality, we get that given $\delta > 0$ for every $\theta$ there are a neighborhood $V(\theta)$ and a constant $J$ such that $d(\varphi(\cdot, \cdot, \theta), \varphi^j(\cdot, \cdot, \theta')) < \delta$ all $\theta' \in V(\theta)$ and all $j \geq J$. Moreover, for every function $f$ and $\epsilon > 0$ this neighborhood $V(\theta)$ and $\delta > 0$ can be chosen small enough so that $|\int f(s)\mu_\theta(ds) - \int f(s)\mu_{\theta'}^j(ds)| < \epsilon$ for all $j \geq J$ and all $\theta' \in V(\theta)$, see [30], Corollary 1. It is now easy to see that the sequence $\{G(\int f(s)\mu_\theta^j(ds), \bar{f})\}_{j \geq 1}$ converges uniformly to $G(\int f(s)\mu_\theta(ds), \bar{f})$. Therefore, the corresponding sequence of minimizers $\{\theta^j\}_{j \geq 1}$ must converge to $\theta^0$. □

PROOF OF PROPOSITION 6.1. This equation for the value of the stock market is obtained from standard arguments using the no-arbitrage conditions along an equilibrium path (e.g., see [18]). □

**Acknowledgment.** This paper was originally written during a long visit of the author to the Departamento de Economia of Universidad Carlos III de Madrid in year 2005.

Department of Economics
University of Miami
5250 University Drive, 517 Jenkings Building
Coral Gables, Florida 33124
USA
E-mail: m.santos2@miami.edu